\documentclass[a4paper]{amsart}
\usepackage{graphicx}
\usepackage{amsmath}
\usepackage{amssymb}
\usepackage{amsfonts}
\usepackage{amsthm}
\usepackage{eucal}
\usepackage{tikz}
\usetikzlibrary{trees, positioning, calc, cd}

\theoremstyle{plain}
\newtheorem{theorem}{Theorem}[section]

\newtheorem*{proposition*}{Proposition}
\newtheorem*{corollary*}{Corollary}

\theoremstyle{definition}

\newtheorem{example}[theorem]{Example}

\makeatletter
\@addtoreset{theorem}{section}
\@addtoreset{subsection}{section}
\makeatother

\begin{document}

\title{Partition complexes and trees}

\author{Gijs Heuts}
\author{Ieke Moerdijk}

\date{}

\begin{abstract}
We construct a homotopy initial functor from the partition complex of a finite set $A$ to a category of trees with leaves labelled by $A$. As an application, this provides an equivalence between different bar constructions of an operad. In the differential graded case, this gives a very elementary proof of an equivalence originally due to Fresse.
\end{abstract}

\maketitle

\section{Introduction}

Let $A$ be a finite set. The partitions of $A$ form a partially ordered set under refinement. This poset has an initial object (the indiscrete partition) and a final object (the discrete partition); removing these two extreme partitions leaves a poset for which we write $\mathbf{P}(A)$. The nerve $N\mathbf{P}(A)$ is a simplicial set called the \emph{partition complex} of $A$. These partition complexes play a central role in the study of bar-cobar (or Koszul) duality for operads \cite{fresse2003koszul,ching}, as we will review in Section \ref{sec:bar}.


In this short note we give a precise and elementary relation between the partition complex $N\mathbf{P}(A)$ and a certain category $\mathbf{T}(A)$ of rooted trees with set of leaves $A$. This category has a morphism $f\colon T \rightarrow S$ if the tree $T$ can be obtained from $S$ by contracting a set of inner edges of $S$. The category $\mathbf{T}(A)$ has an initial object, namely the `corolla' $C_A$ with leaves $A$ and no internal edges. Removing this object from $\mathbf{T}(A)$ gives a subcategory that we denote $\mathbf{T}(A)^+$.

We write $\mathbf{\Delta}/N\mathbf{P}(A)$ for the category of simplices of the simplicial set $N\mathbf{P}(A)$, which can be thought of as a category of layered trees. It has the same weak homotopy type as $N\mathbf{P}(A)$ itself (see Section \ref{sec:partitioncomplexes}), and our result is as follows:

\begin{theorem}
\label{thm:mainthm}
There is a functor
\begin{equation*}
\varphi\colon \mathbf{\Delta}/N\mathbf{P}(A) \rightarrow \mathbf{T}^+(A)
\end{equation*}
with the property that for every object $T$ of $\mathbf{T}^+(A)$, the slice category $\varphi/T$ is weakly contractible. In particular, $\varphi$ induces a weak homotopy equivalence between the classifying spaces of these two categories.
\end{theorem} 

The second statement of the theorem is a consequence of the first by Quillen's Theorem A \cite{quillen1973higher}. However, the first statement is significantly sharper than just $\varphi$ inducing a homotopy equivalence of classifying spaces. A functor $\varphi\colon \mathbf{D} \rightarrow \mathbf{C}$ with the property that $\varphi/c$ is a weakly contractible category for each object $c$ of $\mathbf{C}$ is called \emph{homotopy initial}. The terminology refers to the fact that precomposing by such a functor $\varphi$ preserves homotopy limits. 

The essential content of this paper consists of Sections \ref{sec:partitioncomplexes} and \ref{sec:mainthm}, proving Theorem \ref{thm:mainthm}. We will see that our result implies Fresse's comparison of different versions of the bar construction \cite{fresse2003koszul,livernet2012koszul} for a differential graded operad (see also \cite{vallette2007homology}). In Section \ref{sec:bar} we illustrate some of the uses of the category $\mathbf{T}(A)$ of trees and in particular see how it is related to Ching's description of the bar construction of a topological operad \cite{ching}. In the appendix we explain some of the properties of homotopy initial functors.

\section{Partition complexes and layered trees}
\label{sec:partitioncomplexes}

Fix a finite set $A$. We introduced its partition complex $N\mathbf{P}(A)$ above. The $p$-simplices of this simplicial set are chains of length $p$ of nontrivial partitions, which can be thought of as \emph{layered trees} with $p$ internal layers. For example, the partitions
\begin{equation*}
(abcde)(f) \leq (ab)(cde)(f)
\end{equation*}
of the set $A = \{a,b,c,d,e,f\}$ define a $1$-simplex in $N\mathbf{P}(A)$ represented by the layered tree
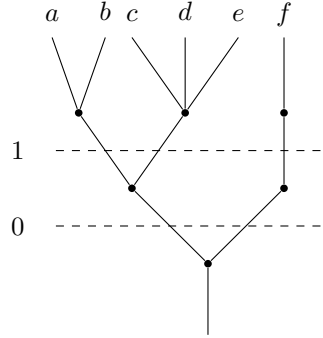
\[
\begin{tikzpicture} 
[level distance=10mm, 
every node/.style={fill, circle, minimum size=.1cm, inner sep=0pt}, 
level 1/.style={sibling distance=20mm}, 
level 2/.style={sibling distance=20mm}, 
level 3/.style={sibling distance=14mm},
level 4/.style={sibling distance=7mm}]

\node (tree)[style={color=white}] {} [grow'=up] 
child {node (level1) {} 
	child{ node (level2) {}
		child{ node (level3) {}
			child
			child
		}
		child{ node {}
			child
			child
			child
		}
	}
	child{ node {}
		child{ node {}
			child
		}
	}
};

\tikzstyle{every node}=[]

\draw[dashed] ($(level1) + (-2cm, .5cm)$) -- ($(level1) + (1.5cm, .5cm)$);
\draw[dashed] ($(level1) + (-2cm, 1.5cm)$) -- ($(level1) + (1.5cm, 1.5cm)$);

\node at ($(level1) + (-2.5cm, 1.5cm)$) {$1$};
\node at ($(level1) + (-2.5cm, .5cm)$) {$0$};
\node at ($(level3) + (-.35cm, 1.3cm)$){$a$};
\node at ($(level3) + (.35cm, 1.35cm)$){$b$};
\node at ($(level3) + (.7cm, 1.3cm)$){$c$};
\node at ($(level3) + (1.4cm, 1.35cm)$){$d$};
\node at ($(level3) + (2.1cm, 1.3cm)$){$e$};
\node at ($(level3) + (2.7cm, 1.3cm)$){$f$};

\end{tikzpicture} 
\]
Here the layer labelled 1 represents the partition $(ab)(cde)(f)$, whereas the layer labelled 0 represents $(abcde)(f)$. Let us fix some terminology.


(1) A simplex in $N\mathbf{P}(A)$ is nondegenerate if it has no layer at which all vertices are unary (meaning they have only one incoming edge, if we orient our trees towards the root). A simplex is called \emph{elementary} if it is has exactly one non-unary vertex at each layer. As an example, observe that the picture above does not represent an elementary simplex, but the trees
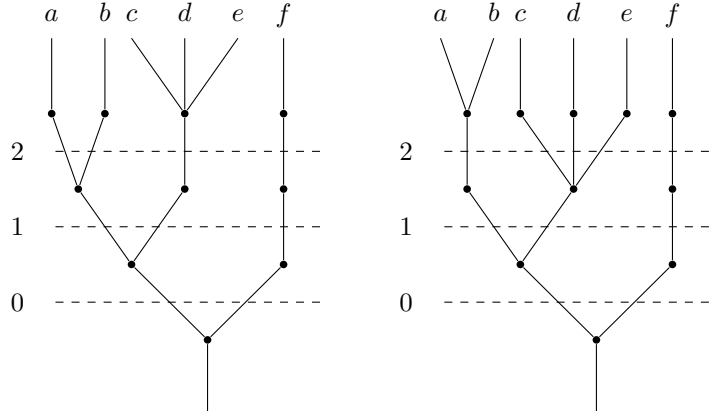
\[
\begin{tikzpicture} 
[level distance=10mm, 
every node/.style={fill, circle, minimum size=.1cm, inner sep=0pt}, 
level 1/.style={sibling distance=20mm}, 
level 2/.style={sibling distance=20mm}, 
level 3/.style={sibling distance=14mm},
level 4/.style={sibling distance=7mm}]

\node (lefttree)[style={color=white}] {} [grow'=up] 
child {node (level1) {} 
	child{ node (level2) {}
		child{ node (level3) {}
			child{ node {}
				child
			}
			child{ node {}
				 child
			}
		}
		child{ node {}
			child{ node {}
				child
				child
				child
			}
		}
	}
	child{ node {}
		child{ node {}
			child{ node {}
				child
			}
		}
	}
};

\node (righttree)[style={color=white}, right = 5cm of lefttree] {} [grow'=up] 
child {node (blevel1) {} 
	child{ node (blevel2) {}
		child{ node (blevel3) {}
			child{ node {}
				child
				child
			}
		}
		child{ node {}
			child{ node {}
				child
			}
			child{ node {}
				child
			}
			child{ node {}
				child
			}
		}
	}
	child{ node {}
		child{ node {}
			child{ node {}
				child
			}
		}
	}
};

\tikzstyle{every node}=[]

\draw[dashed] ($(level1) + (-2cm, .5cm)$) -- ($(level1) + (1.5cm, .5cm)$);
\draw[dashed] ($(level1) + (-2cm, 1.5cm)$) -- ($(level1) + (1.5cm, 1.5cm)$);
\draw[dashed] ($(level1) + (-2cm, 2.5cm)$) -- ($(level1) + (1.5cm, 2.5cm)$);
\draw[dashed] ($(blevel1) + (-2cm, .5cm)$) -- ($(blevel1) + (1.5cm, .5cm)$);
\draw[dashed] ($(blevel1) + (-2cm, 1.5cm)$) -- ($(blevel1) + (1.5cm, 1.5cm)$);
\draw[dashed] ($(blevel1) + (-2cm, 2.5cm)$) -- ($(blevel1) + (1.5cm, 2.5cm)$);

\node at ($(level1) + (-2.5cm, 2.5cm)$) {$2$};
\node at ($(level1) + (-2.5cm, 1.5cm)$) {$1$};
\node at ($(level1) + (-2.5cm, .5cm)$) {$0$};
\node at ($(level3) + (-.35cm, 2.3cm)$){$a$};
\node at ($(level3) + (.35cm, 2.35cm)$){$b$};
\node at ($(level3) + (.7cm, 2.3cm)$){$c$};
\node at ($(level3) + (1.4cm, 2.35cm)$){$d$};
\node at ($(level3) + (2.1cm, 2.3cm)$){$e$};
\node at ($(level3) + (2.7cm, 2.3cm)$){$f$};

\node at ($(blevel1) + (-2.5cm, 2.5cm)$) {$2$};
\node at ($(blevel1) + (-2.5cm, 1.5cm)$) {$1$};
\node at ($(blevel1) + (-2.5cm, .5cm)$) {$0$};
\node at ($(blevel3) + (-.35cm, 2.3cm)$){$a$};
\node at ($(blevel3) + (.35cm, 2.35cm)$){$b$};
\node at ($(blevel3) + (.7cm, 2.3cm)$){$c$};
\node at ($(blevel3) + (1.4cm, 2.35cm)$){$d$};
\node at ($(blevel3) + (2.1cm, 2.3cm)$){$e$};
\node at ($(blevel3) + (2.7cm, 2.3cm)$){$f$};

\end{tikzpicture} 
\]
do represent elementary simplices, namely the two chains of partitions
\begin{eqnarray*}
(abcde)(f) \leq (ab)(cde)(f) \leq (a)(b)(cde)(f), \\
(abcde)(f) \leq (ab)(cde)(f) \leq (ab)(c)(d)(e)(f).
\end{eqnarray*}
As this example illustrates, any nondegenerate simplex in $N\mathbf{P}(A)$ is a face of an elementary one, but the latter need not be unique. 

(2) The category of simplices of $N\mathbf{P}(A)$ is denoted $\mathbf{\Delta}/N\mathbf{P}(A)$. Objects of this category are \emph{layered trees} (with $A$ as set of leaves), morphisms are compositions of face operations deleting a layer and contracting all edges in that layer, and degeneracy operations inserting a new layer consisting entirely of unary vertices. There is a functor
\begin{equation*}
\zeta\colon \mathbf{\Delta}/N\mathbf{P}(A) \longrightarrow \mathbf{P}(A)
\end{equation*}
sending a chain of partitions to its final element. This functor induces a homotopy equivalence of classifying spaces. Indeed, as is well known, this is the case for any small category $\mathbf{C}$ in place of $\mathbf{P}(A)$, since the projection 
\begin{equation*}
\zeta\colon \mathbf{\Delta}/N\mathbf{C} \longrightarrow \mathbf{C}
\end{equation*}
sending $(c_0 \rightarrow \cdots \rightarrow c_p)$ to $c_p$ satisfies the hypothesis of Quillen's Theorem A again, see \cite{thomason1979homotopy}.

\section{Trees without layers and the main theorem}
\label{sec:mainthm}

There is a more basic category of trees $\mathbf{T}(A)$ whose objects are isomorphism classes of finited rooted trees with $A$ as set of leaves, in which every vertex has at least two incoming edges. Morphisms in this category are compositions of `face maps' creating an inner edge, as in
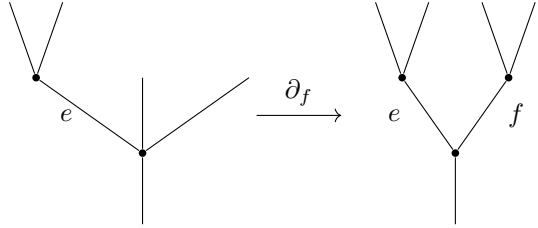
\[
\begin{tikzpicture}
[level distance=10mm, 
every node/.style={fill, circle, minimum size=.1cm, inner sep=0pt}, 
level 1/.style={sibling distance=20mm}, 
level 2/.style={sibling distance=14mm}, 
level 3/.style={sibling distance=7mm},
level 4/.style={sibling distance=7mm}]

\node (lefttree)[style={color=white}] {} [grow'=up] 
child {node (level1) {} 
	child{ node {}
		child
		child
	}
	child
	child
};

\node (righttree)[style={color=white}, right = 4cm of lefttree] {} [grow'=up] 
child {node (blevel1) {} 
	child{ node (blevel2) {}
		child
		child
	}
	child{ node {}
		child
		child
	}
};

\tikzstyle{every node}=[]

\draw[->] ($(level1) + (1.5cm,.5cm)$) -- node[above]{$\partial_f$} ($(blevel1) + (-1.5cm, .5cm)$);

\node at ($(level1) + (-1cm, .5cm)$) {$e$};
\node at ($(blevel1) + (-.8cm, .5cm)$){$e$};
\node at ($(blevel1) + (.8cm, .5cm)$){$f$};

\end{tikzpicture} 
\]
Each such morphism $S \rightarrow T$ maps edges to edges, while it preserves the root and induces the identity on the set of leaves. Note that a morphism $S \rightarrow T$ exists only if $S$ can be obtained from $T$ by contracting a sequence of inner edges and is unique in that case. Thus $\mathbf{T}(A)$ is in fact a poset. This category of trees is the same as the one considered by Ching \cite{ching} and by Hoffbeck and the second author \cite{hoffbeck2021homology}. 

The category $\mathbf{T}(A)$ has an initial object, namely the minimal tree with leaves $A$ connected by a single vertex to its root: we label this tree $C_A$ and refer to it as the corolla with $A$ leaves. We shall write $\mathbf{T}^+(A)$ for the full subcategory of $\mathbf{T}(A)$ obtained by omitting this initial object; in other words, the full subcategory of trees with at least one inner edge.

Forgetting the layers and deleting unary vertices now defines a functor from layered trees to trees,
\begin{equation*}
\varphi\colon \mathbf{\Delta}/N\mathbf{P}(A) \rightarrow \mathbf{T}^+(A),
\end{equation*}
which is the functor of Theorem \ref{thm:mainthm}. The remainder of this section is devoted to its proof.

\begin{proof}[Proof of Theorem \ref{thm:mainthm}]
For an object $T$ of $\mathbf{T}^+(A)$, a \emph{layering} of $T$ is a $p$-simplex $\sigma$ of $N\mathbf{P}(A)$ for which $\varphi(\sigma) = T$. An \emph{elementary} such layering corresponds to a linear ordering of the vertices of $T$ compatible with the partial ordering defined by the tree structure. In particular, $p  = |V(T)| - 2$ in this case, where $V(T)$ is the set of vertices of $T$. If $\tau$ is a face of a layering of $T$, then there exists a map $\varphi(\tau) \rightarrow T$ in $\mathbf{T}^+(A)$. Vice versa, if $\sigma$ is a level tree for which $\varphi(\sigma)$ maps to $T$, then $\sigma$ is a face of a layering of $T$.

Now fix a tree $T$ in $\mathbf{T}^+(A)$. The slice $\varphi/T$ is the category of simplices of a simplicial set $L(T)$ generated by the $(|V(T)|-2)$-simplices representing elementary layerings of $T$. The corolla $C_A$ is not an object of $\mathbf{T}^+(A)$, but we introduce the convention that $L(C_A) = \varnothing$. Also, in the course of the argument below, we will consider trees in $\mathbf{T}^+(B)$ for sets $B$ different from $A$.

We call a vertex $v$ of $T$ a \emph{leaf vertex} if all its incoming edges are leaves of the tree $T$. For such a $v$, let $L^v(T) \subseteq L(T)$ be the simplicial subset generated by simplices which represent elementary layerings of $T$ with $v$ on top. Then
\begin{equation*}
L(T) = \bigcup_{v} L^v(T)
\end{equation*}
with $v$ ranging over leaf vertices of $T$. Moreover, for any such $v$ the simplicial set $L^v(T)$ is clearly a cone on the simplicial set $L(\partial_v T)$, with $\partial_v T$ the tree obtained from $T$ by removing the vertex $v$ and its incoming edges. To be precise, $L^v(T)$ may be identified with the join $L(\partial_v T) \star \Delta[0]$. (This identification results from the observation that any layering of $\partial_v T$, or a face of it, may be naturally upgraded to an object of $L^v(T)$ by attaching the corolla with vertex $v$ at the top.) Similarly, if $v_1, \ldots, v_k$ is a finite set of leaf vertices of $T$, then the intersection
\begin{equation*}
L^{v_1} (T) \cap \cdots \cap L^{v_n}(T) 
\end{equation*}
is a cone on $L(\partial_{v_1, \ldots, v_n} T)$, where $\partial_{v_1, \ldots, v_n} T$ is obtained from $T$ by removing all of the vertices $v_1, \ldots, v_n$ and their incoming edges. Thus, the $L^v(T)$ together form a cover of $L(T)$ by weakly contractible simplicial sets, all of whose finite intersections are again contractible. So $L(T)$ is itself weakly contractible, which proves the theorem.
\end{proof}

\begin{example}
\label{ex:labelledpartitions}
An operad $\mathbf{O}$ in the category of sets defines a presheaf $N\mathbf{O}$ on $\mathbf{T}^+(A)$ called its \emph{nerve}. This presheaf assigns to any tree $T$ the set of labellings of the vertices of $T$ by operations in $\mathbf{O}$ of arity corresponding to the number of inputs of the vertex. The functoriality with respect to edge contractions uses the composition of operations in $\mathbf{O}$. The category of elements of $N\mathbf{O}$ is denoted
\begin{equation*}
\mathbf{T}^+(A,\mathbf{O}) := \mathbf{T}^+(A)/N\mathbf{O}.
\end{equation*}
The presheaf $N\mathbf{O}$ pulls back to a presheaf $\varphi^*(N\mathbf{O})$ on $\mathbf{\Delta}/N\mathbf{P}(A)$ whose category of elements we denote $\mathbf{\Delta}/N\mathbf{P}(A)_{\mathbf{O}}$. It is the category of elements of a simplicial set $N\mathbf{P}(A)_{\mathbf{O}}$. This simplicial set is a \emph{labelled partition complex} of the kind considered in \cite{fresse2003koszul}. As explained in the appendix, our theorem implies that the functor
\begin{equation*}
\mathbf{\Delta}/N\mathbf{P}(A)_{\mathbf{O}} \longrightarrow \mathbf{T}^+(A,\mathbf{O})
\end{equation*}
is again homotopy initial, so in particular a weak homotopy equivalence. The same applies if $\mathbf{O}$ is an operad in simplicial sets, or in chain complexes, for example. (As explained in the appendix, really any symmetric monoidal model category or $\infty$-category could be used instead.) This gives a different and more direct proof of the quasi-isomorphism of Fresse between different bar constructions for an operad proved in \cite{fresse2003koszul} and mentioned in the introduction, see also the discussion in the next section. Note that Fresse's proof goes via a construction in the opposite direction, from trees to layered trees, which he calls `levelization'.
\end{example}

\section{More on the category of trees}
\label{sec:bar}

In this concluding section we collect some remarks on the category of trees $\mathbf{T}(A)$ described above, which should illustrate its relevance and some of its basic properties.

(1) As $A$ ranges over finite sets, the categories $\mathbf{T}(A)$ (as well as $\mathbf{T}^+(A)$) form an operad in the category of small categories. Indeed, for $a\in A$, a composition map of the form
\begin{equation*}
\mathbf{T}(B) \circ_a \mathbf{T}(A) \rightarrow \mathbf{T}(A \circ_a B)
\end{equation*}
is given by grafting trees with leaves $B$ onto the leaf $a$ of a tree with leaves $A$. Here $A \circ_a B$ denotes the set $(A-\{a\}) \amalg B$. An algebra for this operad is a category with a colax symmetric monoidal structure, as has been emphasized by Hinich.

(2) The homotopy type of the classifying space of $\mathbf{T}^+(A)$ (or equivalently of the partition complex $N\mathbf{P}(A)$) is a wedge of $(|A|-1)!$ spheres of dimension $|A|-3$ (cf. \cite{robinson1996tree}).

(3) Let $\mathcal{E}$ be a symmetric monoidal simplicial model category (or symmetric monoidal $\infty$-category) with a zero object. Then there is a `trivial' operad $\mathbf{1}$ having the monoidal unit in arity 1 and the zero object in all others. If $\mathbf{O} \rightarrow \mathbf{1}$ is an augmented operad in $\mathcal{E}$ one may form its simplicial bar construction, with object of $p$-simplices
\begin{equation*}
B(\mathbf{1}, \mathbf{O}, \mathbf{1})_p = \mathbf{O}\circ \cdots \circ \mathbf{O}.
\end{equation*}
Here there are $p$ copies of $\mathbf{O}$ on the right, the symbol $\circ$ denotes the usual composition product of symmetric sequences, the inner face maps use the operad structure of $\mathbf{O}$, and the outer face maps apply the augmentation. The realization of this simplicial object in $\mathcal{E}$ is (one possible version of) the bar construction $B\mathbf{O}$ of the operad $\mathbf{O}$.

It is not difficult to see (cf. \cite[Section 4.3]{fresse2003koszul}) that this $B\mathbf{O}$ can be described explicitly in terms of the labelled partition complexes introduced above. Indeed, the term $\mathbf{O}^{\circ p}(A) = B(\mathbf{1}, \mathbf{O}, \mathbf{1})_p(A)$ consists of layered trees of height $p$ with leaves $A$, and each internal vertex labelled by an operation of $\mathbf{O}$. Combining this with our result, one deduces fairly easily that $B\mathbf{O}(A)$ is a suspension of the cofiber of the map
\begin{equation*}
\mathrm{hocolim}_{T \in \mathbf{T}^+(A)} N\mathbf{O}(T) \rightarrow \mathrm{hocolim}_{T \in \mathbf{T}(A)} N\mathbf{O}(T) \simeq \mathbf{O}(A).
\end{equation*}
This gives an alternative description for the simplicial bar construction in terms of trees. In the case where $\mathcal{E}$ is the category of chain complexes over a commutative ring, this reproduces the comparison of bar constructions given by Fresse \cite{fresse2003koszul}. If $\mathcal{E}$ is the category of topological spaces or spectra, Ching \cite{ching} observes that the simplicial bar construction is even \emph{homeomorphic} to the one defined in terms of trees; our result is a bit weaker, but applies more generally.

(4) If again $\mathbf{O}$ is an operad in a suitable simplicial model category, then the term $\mathrm{hocolim}_{T \in \mathbf{T}(A)} N\mathbf{O}$ is essentially the Boardman--Vogt $W$-construction $W\mathbf{O}(A)$ (as described in \cite{berger2006boardman}). Indeed, the latter can be represented by `the space' of trees with leaves $A$, vertices labelled by operations of $\mathbf{O}$ and inner edges labelled by a length $t$ in the 1-simplex $\Delta^1$. If such a length is 0, then the resulting tree is identified with the one where that inner edge is contracted and the operations corresponding to its endpoints are composed. Viewed in this way, the subobject $\mathrm{hocolim}_{T \in \mathbf{T}^+(A)} N\mathbf{O}(T)$ can be thought of as the \emph{decomposable} operations in $W\mathbf{O}(A)$, i.e., the ones arising as a composition of two or more operations in the operad $W\mathbf{O}$. Writing $\mathrm{Indec}(W\mathbf{O}(A))$ for the quotient of $W\mathbf{O}(A)$ by the decomposable operations, the previous remark implies an equivalence of symmetric sequences
\begin{equation*}
B\mathbf{O} \simeq \Sigma\mathrm{Indec}(W\mathbf{O}).
\end{equation*}
In the context of topological operads, this observation goes back to Salvatore.

(5) The Fulton--MacPherson operad $\mathbf{FM}_d$ has as its terms $\mathbf{FM}_d(A)$ certain manifolds with corners that compactify the configuration spaces $\mathrm{Conf}_d(A)$ of $|A|$ points in $\mathbb{R}^d$, or rather their quotients by translation and dilation (cf. the discussion in \cite[Section 1.8]{heutsmoerdijk}). These manifolds are naturally stratified, with strata indexed by the objects of the category $\mathbf{T}(A)$. The stratum corresponding to a tree is homeomorphic to a product over the vertices of $T$ of the  configuration spaces corresponding to the inputs of each vertex; the corolla $C_A$ corresponds to the maximal open stratum. The operad $\mathbf{FM}_d$ is weakly equivalent to the operad $\mathbf{E}_d$ of little $d$-cubes but has more favourable cofibrancy properties. 

Using the description of $\mathbf{FM}_d(A)$ as a stratified space, it is not difficult to see that the homotopy colimit 
\begin{equation*}
\mathrm{hocolim}_{\mathbf{T}^+(A)} N\mathbf{FM}_d(-)
\end{equation*}
is equivalent to the boundary of the manifold $\mathbf{FM}_d(A)$. Thus, one deduces from (3) above that the terms of the bar construction (and hence of the bar construction of the little $d$-disks operad $\mathbf{E}_d$) are equivalent to the spaces
\begin{equation*}
\Sigma\bigl(\mathbf{FM}_d(A)/\partial\mathbf{FM}_d(A)\bigr).
\end{equation*}
The reader might compare this to the remark immediately following \cite[Lemma 9.5]{ching2020koszul}. This reference contains a much more detailed analysis of the bar construction of the Fulton--MacPherson operad and includes a proof of the `self-duality' of the $\mathbf{E}_d$-operad.

\section*{Appendix}
\label{appendix}

This appendix is a service to the reader not familiar with homotopy initial functors, illustrating some useful properties of this notion. None of this material is original. A standard account in the context of model categories is \cite[Section 19.6]{hirschhorn} or \cite[Section 4.1]{htt} in the context of $\infty$-categories.

If $\varphi\colon \mathbf{D} \rightarrow \mathbf{C}$ is a functor and $c$ an object of $\mathbf{C}$, then $\varphi/c$ (or simply $\mathbf{D}/c$, leaving $\varphi$ implicit) denotes the category with objects the pairs consisting of an object $d$ of $\mathbf{D}$ and a morphism $\alpha\colon\varphi(d) \rightarrow c$, and with morphisms from $\alpha$ to another $\alpha'\colon \varphi(d') \rightarrow c$ the morphisms $\beta\colon d \rightarrow d'$ in $\mathbf{D}$ making the resulting triangle
\[
\begin{tikzcd}
\varphi(d) \ar{dr}[swap]{\alpha}\ar{rr}{\varphi(\beta)} && \varphi(d') \ar{dl}{\alpha'} \\
& c &
\end{tikzcd}
\]
commute. The functor $\varphi$ is called \emph{homotopy initial} (or \emph{homotopy left cofinal} by some) if for each object $c$ of $\mathbf{C}$, the category $\varphi/c$ has contractible classifying space.

Quillen's Theorem A \cite{quillen1973higher} states that if $\varphi$ is homotopy initial, then it induces a homotopy equivalence of classifying spaces $B\mathbf{D} \simeq B\mathbf{C}$. However, the properties of a homotopy initial functor are much stronger than this. As suggested by the terminology, restriction along $\varphi$ preserves arbitrary homotopy limits, as explained in (e) below. (Alternatively, restriction along $\varphi\colon \mathbf{D}^{\mathrm{op}} \rightarrow \mathbf{C}^{\mathrm{op}}$ preserves homotopy colimits.) Let us list several useful properties of such functors, most of which can be seen as special instances of the last example (e).

(a) If $X$ is a presheaf of sets on $\mathbf{C}$, with its restriction $\varphi^*X$ a presheaf on $\mathbf{D}$, then the induced functor
\begin{equation*}
\mathbf{D}/\varphi^*X \rightarrow \mathbf{C}/X
\end{equation*}
between the respective categories of elements is again homotopy initial (as is easily verified), hence induces a homotopy equivalence
\begin{equation*}
B(\mathbf{D}/\varphi^*X) \xrightarrow{\simeq} B(\mathbf{C}/X).
\end{equation*}
Note that this map is a pullback of the map $B\mathbf{D} \rightarrow B\mathbf{C}$, but that the projection of the codomain to $B\mathbf{C}$ is rarely a fibration.

(b) The same applies more generally for a fibered category $p\colon \mathbf{E} \rightarrow \mathbf{C}$ and its pullback $\varphi^* \mathbf{E} \rightarrow \mathbf{D}$. To prove this, one shows that for $e$ an object of $\mathbf{E}$, the classifying space of $\varphi^*\mathbf{E}/e$ is weakly equivalent to that of $\mathbf{D}/p(e)$ by another application of Quillen's Theorem A. (The directions are important here; the result would not generally hold for opfibered categories.)

(c) Generalizing (a), one could take $X_\bullet$ to be a simplicial presheaf on $\mathbf{C}$. Then
\begin{equation*}
\mathbf{D}/\varphi^*X_\bullet \rightarrow \mathbf{C}/X_\bullet
\end{equation*}
is a map of bisimplicial sets that is a weak equivalence in each fixed simplicial degree (of $X_\bullet$), hence induces a weak equivalence on diagonals. 

(d) If $A$ is a presheaf of abelian groups on $\mathbf{C}$, then $\varphi$ induces an isomorphism $H_*(\mathbf{C};A) \rightarrow H_*(\mathbf{D};\varphi^*A)$ in homology. (Here homology is defined as the collection of left derived functors of the colimit functor from abelian presheaves on $\mathbf{C}$ to the category of abelian groups.) For this result it is not necessary to require that $A$ is locally constant, as for ordinary weak equivalences. Indeed, if $P_\bullet \rightarrow A$ is a resolution of $A$ where $P_\bullet = \mathbb{Z}[X_\bullet]$ is free on a simplicial presheaf $X_\bullet$, then the presheaf $\varphi^*(P_\bullet) = \mathbb{Z}[\varphi^*X_\bullet]$ is a similar resolution of $X_\bullet$. The result then follows because for each simplicial degree $p$, the map $H_*(\mathbf{C};\mathbb{Z}[X_p]) \rightarrow H_*(\mathbf{D};\mathbb{Z}[\varphi^*X_p])$ can be identified with
\begin{equation*}
H_*(\mathbf{C}/X_p; \mathbb{Z}\Bigr) \xrightarrow{\simeq} H_*(\mathbf{D}/\varphi^*(X_p);\mathbb{Z}),
\end{equation*}
which is an isomorphism by case (a).

(e) Consider the category $\mathcal{E}^{\mathbf{C}^{\mathrm{op}}}$ of presheaves with values in a simplicial model category $\mathcal{E}$, and similarly $\mathcal{E}^{\mathbf{D}^{\mathrm{op}}}$. These presheaf categories can be equipped with the projective model structures and $\varphi$ induces a Quillen adjunction
\[
\begin{tikzcd}
\mathcal{E}^{\mathbf{D}^{\mathrm{op}}} \ar[shift left]{r}{\varphi_!} & \mathcal{E}^{\mathbf{C}^{\mathrm{op}}}. \ar[shift left]{l}{\varphi^*}
\end{tikzcd}
\]
Constructing a similar adjunction for the trivial functor $\mathbf{C} \rightarrow 1$ gives the Quillen adjunction
\[
\begin{tikzcd}
\mathcal{E}^{\mathbf{C}^{\mathrm{op}}} \ar[shift left]{r}{\mathrm{colim}_{\mathbf{C}}} & \mathcal{E} \ar[shift left]{l}{\mathrm{const}}
\end{tikzcd}
\]
and similarly for $\mathbf{D}$. If $\varphi$ is homotopy initial, then there is a natural isomorphism of derived functors (cf. \cite[Section 19.6]{hirschhorn})
\begin{equation*}
\mathbf{L}\mathrm{colim}_{\mathbf{D}} \circ \mathbf{R}\varphi^* \cong \mathbf{L}\mathrm{colim}_{\mathbf{C}}.
\end{equation*}
In other words, the restriction functor $\mathbf{R}\varphi^*$ preserves homotopy colimits.

The proof of this is rather straightforward from the usual Bousfield--Kan formula for homotopy colimits, which states that for an $\mathcal{E}$-valued presheaf $X$ on $\mathbf{C}$ taking values in cofibrant objects, its homotopy colimit may be computed as a tensor product
\begin{equation*}
X(-) \otimes_{\mathbf{C}} N(\mathbf{C}/-).
\end{equation*}
Similarly, one computes the homotopy colimit of $\varphi^*X$ as
\begin{equation*}
X(\varphi(-)) \otimes_{\mathbf{D}} N(\mathbf{D}/-).
\end{equation*}
The latter is isomorphic to
\begin{equation*}
X(-) \otimes_{\mathbf{C}} N(\mathbf{C}/\varphi(-)),
\end{equation*}
and the natural comparison map $N(\mathbf{C}/\varphi(-)) \rightarrow N(\mathbf{C}/-)$ is a weak homotopy equivalence of simplicial diagrams on $\mathbf{C}$ by the assumption that $\varphi$ is homotopy initial. This implies the result.

Examples (a) and (d) can be seen as special cases of this preservation of homotopy colimits: for example (a), this relies on the fact that the classifying space of the category of elements $\mathbf{C}/X$ is a model for the homotopy colimit of the functor
\begin{equation*}
\mathbf{C}^{\mathrm{op}} \rightarrow \mathbf{Sets} \xrightarrow{\mathrm{const}} \mathbf{sSets}.
\end{equation*}

(f) There is an $\infty$-categorical version of (e) above, stating that precomposition by $\varphi$ preserves colimits with values in a fixed $\infty$-category $\mathcal{E}$ (cf. \cite[Theorem 4.1.3.1]{htt}). Also, example (b) can be seen as a special instance of the fact the pullback of a homotopy initial functor along a cartesian fibration is again homotopy initial \cite[Proposition 4.1.2.15]{htt}.

\bibliographystyle{plain}
\bibliography{biblio}

\end{document}